\documentclass[12pt,reqno]{amsart}
\newtheorem{thm}{Theorem}[section]
\newtheorem{lma}[thm]{Lemma}

\newtheorem{cor}[thm]{Corollary}

\newtheorem{rem}[thm]{Remark}

\newcommand{\R}{{\mathbb R}}

\newcommand{\N}{{\mathbb N}}

\newcommand{\cF}{{\mathcal F}}

\newcommand{\eps}{{\varepsilon}}
\newcommand{\pf}{{\noindent \bf Proof. }}

\newcommand{\dist}{{\rm dist\ }}

\newcommand{\supp}{{\rm supp\ }}
\newcommand\dd{{\,\mathrm d}}

\newcommand\vol{\operatorname{vol}}

\numberwithin{equation}{section}
\begin{document}

\title[Optimal stability estimate]{Optimal stability estimate of the inverse boundary value problem by partial measurements}

\author[Heck]{Horst Heck}
\address{Technische Universit\"at Darmstadt, FB Mathematik, AG 4, Schlossgartenstr. 7, D-64289 Darmstadt, Germany}
\email{heck@mathematik.tu-darmstadt.de}
\author[Wang]{Jenn-Nan Wang}
\address{Department of Mathematics, National Taiwan University,
Taipei 106, Taiwan}
\email{jnwang@math.ntu.edu.tw}
\thanks{The second author was supported in part by the National Science Council of Taiwan (NSC 95-2115-M-002-003).}

%
\begin{abstract}
In this work we establish \emph{log} type stability estimates for
the inverse potential and conductivity problems with partial
Dirichlet-to-Neumann map, where the Dirichlet data is homogeneous on
the inaccessible part. This result, to some extent, improves our
former result on the partial data problem \cite{HW06} in which
\emph{log-log} type estimates were derived.
\end{abstract}

\maketitle

\section{Introduction}

In this paper we study the stability question of the inverse
boundary value problem for the Schr\"odinger equation with a
potential and the conductivity equation by partial Cauchy data. This
type of inverse problem with full data, i.e., Dirichlet-to-Neumann
map, were first proposed by Calder\'on \cite{C80}. For three or
higher dimensions, the uniqueness issue was settled by Sylvester and
Uhlmann \cite{SU87} and a reconstruction procedure was given by
Nachman \cite{N88}. For two dimensions, Calder\'on's problem was
solved by Nachman \cite{N96} for $W^{2,p}$ conductivities and by
Astala and P\"aiv\"arinta \cite{AP06} for $L^{\infty}$
conductivities. This inverse problem is known to be ill-posed. A
log-type stability estimate was derived by Alessandrini \cite{A88}.
On the other hand, it was shown by Mandache \cite{M01} that the
log-type estimate is optimal.

All results mentioned above are concerned with the full data.
Recently, the inverse problem with partial data has received lots of
attentions \cite{GU01}, \cite{IU04}, \cite{BU02}, \cite{KSU05},
\cite{FKSU07}, \cite{I07}. A {\it log-log} type stability estimate
for the inverse problem with partial data was derived by the authors
in \cite{HW06}. The method in \cite{HW06} was based on \cite{BU02} and
a stability estimate for the analytic continuation proved in
\cite{Ve99}. We believe that the log type estimate should be the
right estimate for the inverse boundary problem, even with partial
data. In this paper, motivated by the uniqueness proof in Isakov's
work \cite{I07}, we prove a {\it log} type estimate for the inverse
boundary value problem under the same {\it a priori} assumption on
the boundary as given in \cite{I07}. Precisely, the inaccessible
part of the boundary is either a part of a sphere or a plane. Also,
one is able to use zero data on the inaccessible part of the
boundary. The strategy of the proof in \cite{I07} follows the
framework in \cite{SU87} where complex geometrical optics solutions
are key elements. A key observation in \cite{I07} is that when
$\Gamma_0$ is a part of a sphere or a plane, we are able to use a
reflection argument to guarantee that complex geometrical optics
solutions have homogeneous data on $\Gamma_0$.

Let $n\ge 3$ and $\Omega\subset\R^n$ be an open domain with smooth
boundary $\partial\Omega$. Given $q\in L^{\infty}(\Omega)$, we
consider the boundary value problem:
\begin{align}\label{bvp}
\begin{split}
(\Delta-q)u & =0\quad\text{in}\;\Omega\\
u &=f\quad\text{on}\;\partial\Omega,
\end{split}
\end{align}
where $f\in H^{1/2}(\partial\Omega)$. Assume that $0$ is not a
Dirichlet eigenvalue of $\Delta-q$ on $\Omega$. Then \eqref{bvp}
has a unique solution $u\in H^1(\Omega)$. The usual definition of
the Dirichlet-to-Neumann map is given by
\[
\Lambda_qf=\partial_{\nu}u|_{\partial\Omega}
\]
where  $\partial_{\nu}u=\nabla u\cdot\nu$ and $\nu$ is the unit
outer normal of $\partial\Omega$.

Let $\Gamma_0\subset\partial\Omega$ be an open part of the boundary
of $\Omega$. We set $\Gamma=\partial\Omega\setminus\Gamma_0$. We
further set $H_0^{1/2}(\Gamma):=\{f\in H^{1/2}(\partial\Omega):
\supp f \subset \Gamma\}$ and $H^{-1/2}(\Gamma)$ the dual space of
$H^{1/2}_0(\Gamma)$. Then the partial Dirichlet-to-Neumann map
$\Lambda_{q,\Gamma}$ is defined as
\[
\Lambda_{q,\Gamma} f:= \partial_\nu u|_{\Gamma}\in H^{-1/2}(\Gamma)
\]
where $u$ is the unique weak solution of \eqref{bvp} with Dirichlet
Data $f\in H_0^{1/2}(\Gamma)$. In what follows, we denote the
operator norm by
\[
\|\Lambda_{q,\Gamma}\|_{\ast}:=
\|\Lambda_{q,\Gamma}\|_{H^{1/2}_0(\Gamma)\to H^{-1/2}(\Gamma)}
\]

We consider two types of domains in this paper: (a) $\Omega$ is a
bounded domain in $\{x_n<0\}$ and
$\Gamma_0=\partial\Omega\cap\{x_n<0\}$; (b) $\Omega$ is a subdomain of
$B(a,R)$ and $\Gamma_0=\partial B(a,R)\cap\partial\Omega$ with
$\Gamma_0\neq\partial B(a,R)$, where $B(a,R)$ is a ball centered at
$a$ with radius $R$. Denote by $\hat q$ the zero extension of the
function $q$ defined on $\Omega$ to $\R^n$. The main result of the
paper reads as follows:
\begin{thm}\label{main_bdd}
Assume that $\Omega$ is given as in either (a) or (b). Let $N>0$,
$s>\frac{n}{2}$ and $q_j\in H^s(\Omega)$ such that
\begin{equation}\label{apr}
\|q_j\|_{H^s(\Omega)}\le N
\end{equation}
for $j=1,2$, and $0$ is not a Dirichlet eigenvalue of $\Delta-q_j$ for $j=1,2$.
Then there exist constants $C>0$ and $\sigma >0$
such that
\begin{equation}\label{est}
\|q_1-q_2\|_{L^\infty(\Omega)} \leq C \big{|}\log
\|\Lambda_{q_1,\Gamma}-
    \Lambda_{q_2,\Gamma}\|_{*} \big{|}^{-\sigma}
\end{equation}
where $C$ depends on $\Omega, N, n, s$ and
$\sigma$ depends on $n$ and $s$.
\end{thm}

Theorem~\ref{main_bdd} can be generalized to the conductivity
equation. Let $\gamma \in H^s(\Omega)$ with $s>3+\frac n2$ be a
strictly positive function on $\overline{\Omega}$. The equation for
the electrical potential in the interior without sinks or sources is
\begin{align*}
\text{div}(\gamma\nabla u) &=0\quad\text{in}\quad\Omega\\
u &=f\quad\text{on}\quad\partial\Omega.
\end{align*}
As above, we take $f\in H^{1/2}_0(\Gamma)$. The partial
Dirichlet-to-Neumann map defined in this case is
\[
\Lambda_{\gamma,\Gamma}:f\mapsto\gamma\partial_{\nu}u|_{\Gamma}.
\]
\begin{cor}\label{cor1}
Let the domain $\Omega$ satisfy (a) or (b). Assume that
$\gamma_j\ge N^{-1}>0$, $s>\frac n2$, and
\begin{equation}\label{ap}
\|\gamma_j\|_{H^{s+3}(\Omega)}\le N
\end{equation}
for $j=1,2$, and
\begin{equation}\label{sameb}
\partial_{\nu}^{\beta}\gamma_1|_{\Gamma}=\partial_{\nu}^{\beta}\gamma_2|_{\Gamma}\quad\text{on}\quad\partial\Omega,\quad\forall\quad
0\le\beta\le 1.
\end{equation}
Then there exist constants $C>0$ and $\sigma>0$ such that
\begin{equation}\label{est111}
\|\gamma_1-\gamma_2\|_{L^\infty(\Omega)} \leq C \big{|}\log
\|\Lambda_{\gamma_1,\Gamma}-
    \Lambda_{\gamma_2,\Gamma}\|_{*} \big{|}^{-\sigma}
\end{equation}
where $C$ depend on $\Omega, N, n, s$ and $\sigma$ depend on $n,s$.
\end{cor}
\begin{rem}
For the sake of simplicity, we impose the boundary identification
condition \eqref{sameb} on conductivities. However, using the
arguments in \cite{A88} {\rm (}also see \cite{HW06}{\rm )}, this
condition can be removed. The resulting estimate is still in the
form of \eqref{est111} with possible different constant $C$ and
$\sigma$.
\end{rem}

\section{Preliminaries}

We first prove an estimate of the Riemann-Lebesgue lemma for a certain class of functions. Let us define
$$g(y)=\|f(\cdot-y)-f(\cdot)\|_{L^1(\R^n)}$$ for any $f\in L^1(\R^n)$. It is known that
$\lim_{|y|\to 0}g(y)=0.$
\begin{lma}\label{fourier_est}
Assume that $f\in L^1(\R^n)$ and there exist $\delta>0$, $C_0>0$, and $\alpha\in (0,1)$ such that
\begin{equation}\label{holder}
g(y)\le C_0|y|^{\alpha}
\end{equation}
whenever $|y|<\delta$. Then there exists a constant $C>0$ and $\eps_0>0$
 such that for any $0<\varepsilon<\eps_0$ the inequality
\begin{equation}\label{est1}
|\cF f(\xi)| \leq C(\exp(-\pi\varepsilon^2|\xi|^2) +
\varepsilon^\alpha) \end{equation} holds with $C=C(C_0,\|f\|_{L^1},n,\delta,\alpha)$.
\end{lma}
\pf
Let $G(x):= \exp(-\pi|x|^2)$ and set $G_\varepsilon(x):=
\varepsilon^{-n} G(\frac{x}{\varepsilon})$. Then we define
$f_\varepsilon := f\ast G_\varepsilon$. Next we write
\[
|\cF f(\xi)| \leq |\cF f_\varepsilon(\xi)| + |\cF (f_\varepsilon -
f)(\xi)|.
\] 
For the first term on the right hand side we get
\begin{align}\label{est2}
\begin{split}
|\cF f_\varepsilon(\xi)| & \leq |\cF f(\xi)|\cdot |\cF
		G_\varepsilon(\xi)| \\
	& \leq \|f\|_1|\varepsilon^{-n}\varepsilon^n \cF
		G(\varepsilon\xi)| \\
	&\leq \|f\|_1\exp(-\pi \varepsilon^2|\xi|^2).
\end{split}
\end{align}
To estimate the second term, we use the assumption \eqref{holder} and derive
\begin{align*}
|\cF (f_\varepsilon - f)(\xi)| &\leq \|f_\varepsilon - f\|_1\\
    &\leq \int_{\R^n}\int_{\R^n} |f(x-y) - f(x)| G_\varepsilon(y) \dd y
            \dd x\\
    &= \int_{|y|<\delta}\int_{\R^n} |f(x-y) - f(x)| G_\varepsilon(y)\dd x \dd y
            \\
\begin{split}
&+ \int_{|y|\ge\delta}\int_{\R^n} |f(x-y)-f(x)| G_\varepsilon(y)\dd x \dd y
\end{split}\\
    &= I + II.
\end{align*}
In view of \eqref{holder} we can estimate
\begin{align*}
I&=\int_{|y|<\delta}g(y)G_\eps (y) \dd y\\
    &\leq C_0\int_{|y|<\delta}|y|^{\alpha}G_\eps (y)
        \dd y\\
    &= C_0\int_{S^{n-1}}\int_{0}^{\delta} r^{\alpha}\eps^{-n}
        \exp(-\pi \eps^{-2} r^2)r^{n-1} \dd r \dd \psi\\
    &= C_1\int_{0}^{\delta}\eps^{\alpha}u^{\alpha}\eps^{-n}
    	\exp(-u^2)\eps^{n-1}u^{u-1}\eps \dd u\\
    &= C_2\eps^{\alpha}\int_{0}^{\delta}u^{n+\alpha-1}\exp(-u^2) \dd
    	u=C_3\eps^{\alpha},
\end{align*}
where $C_3=C_3(C_0,n,\delta,\alpha)$.

As for II, we obtain that for $\eps$ sufficiently small
\begin{align*}
II&=\int_{|y|\ge\delta}g(y)G_\eps (y)\dd y\\
    &\le 2\|f\|_{L^1}\int_{|y|\ge\delta}G_\eps (y)\dd y\\
    &\leq C_4\|f\|_1 \int_{\delta}^{\infty} \eps^{-n} \exp(-\pi\eps^{-2}
        r^2)r^{n-1} \dd r \\
    &= C_4\|f\|_1\int_{\delta\eps^{-1}}^{\infty} u^{n-1} \exp(-\pi u^2)\dd u \\
    &\leq C_4\|f\|_1 \int_{\delta\eps^{-1}}^{\infty} \exp(-\pi u)\dd u \\
    &\leq C_4\|f\|_1 \frac{1}{\pi} \exp(-\pi\delta \eps^{-1}) \leq C_5
        \eps^{\alpha},
\end{align*}
where $C_5=C_5(\|f\|_{L^1},n,\delta,\alpha)$.
Combining the estimates for $I$, $II$, and \eqref{est2}, we
immediately get \eqref{est1}.
\qed

We now provide a sufficient condition on $f$, defined on $\Omega$,
such that \eqref{holder} in the previous lemma holds.
\begin{lma}\label{lma:holder_est}
Let $\Omega\subset \R^n$ be a bounded domain with $C^1$ boundary.
Let $f\in C^\alpha(\overline{\Omega})$ for some $\alpha\in (0,1)$
and denote by $\hat f$ the zero extension of $f$ to $\R^n$. Then
there exists $\delta>0$ and $C>0$ such that
\[
\|\hat f(\cdot-y)-\hat f(\cdot)\|_{L^1(\R^n)} \leq C|y|^\alpha
\]
for any $y\in \R^n$ with $|y|\leq \delta$.
\end{lma}
\pf
Since $\Omega$ is bounded and of class $C^1$, there exist a finite
number of balls, say $m\in \N$, $B_i(x_i)$ with center $x_i\in
\partial\Omega$, $i=1,\dots,m$ and associated $C^1$-diffeomorphisms
$\varphi_i:B_i(x_i)\to Q$ where $Q=\{x'\in \R^{n-1}: \|x'\|\leq
1\}\times (-1,1)$. Set
$d=\dist(\partial\Omega,\partial(\bigcup_{i=1}^m B_i(x_i)))>0$
and $\tilde\Omega_\eps=\bigcup_{x \in\partial\Omega} B(x,\eps)$,
where $B(x,\eps)$ denotes the ball with center $x$ and radius
$\eps>0$. Obviously, for $\eps<d$, it holds that
$\tilde\Omega_\eps\subset \bigcup_{i=1}^m B_i(x_i)$. Let
$x\in\partial\Omega$ and $0<|y|<\delta\le d$, then for any
$z_1,z_2\in B(x,|y|)\cap B_i(x_i)$ we get that
\[
|\varphi_i(z_1)-\varphi_i(z_2)| \leq
\|\nabla\varphi_i\|_{L^{\infty}} |z_1-z_2| \leq C|y|
\]
for some constant $C>0$. Therefore, $\varphi_i(\tilde\Omega_{|y|}\cap
B_i(x_i)) \subset \{x'\in \R^{n-1}: \|x'\|\leq 1\}\times
(-C|y|,C|y|)$. By the transformation formula this yields
$\vol(\tilde\Omega_{|y|})\leq C|y|$.

Since $|y|<\delta$ we have $\hat f(x-y) -\hat f(x)=0$ for $x
\not \in \Omega\cup \tilde\Omega_{|y|}$. Now we write
\begin{align*}
\| \hat f(\cdot-y) -\hat f\|_{L^1(\R^n)}  =
&\int_{\Omega\setminus
        \tilde\Omega_{|y|}} |\hat f(x-y) -\hat f(x)| \dd x\\
         +&
        \int_{\tilde\Omega_{|y|}} |\hat f(x-y) -\hat f(x)| \dd x\\
    \leq\strut &C\vol(\Omega) |y|^\alpha + 2\|f\|_{L^\infty}
        \vol(\tilde\Omega_{|y|})\\
    \leq\strut &C(|y|^\alpha+|y|) \leq C|y|^\alpha
\end{align*}
for $\delta\leq 1$.
\qed

Now let $q_1$ and $q_2$ be two potentials and their corresponding
partial Dirichlet-to-Neumann maps are denoted by
$\Lambda_{1,\Gamma}$ and $\Lambda_{2,\Gamma}$, respectively. The
following identity plays a key role in the derivation of the stability
estimate.
\begin{lma}\label{lma_id}
Let $v_j$ solve \eqref{bvp} with $q=q_j$ for $j=1,2$. Further assume
that $v_1=v_2=0$ on $\Gamma_0$. Then
\[
\int_\Omega (q_1-q_2)v_1\overline{v_2}\dd x = \langle
(\Lambda_{1,\Gamma}-\Lambda_{2,\Gamma})v_1,v_2 \rangle
\]
\end{lma}
\pf
Let $u_2$ denote the solution of \eqref{bvp} with $q=q_2$ and
$u_2=v_1$ on $\partial\Omega$. Therefore
\begin{align*}
\int_\Omega \nabla v_1\cdot\overline{\nabla v_2} +
q_1v_1\overline{v_2} \dd x &=
\langle \partial_\nu v_1, v_2\rangle\\
\int_\Omega \nabla u_2\cdot\overline{\nabla v_2} +
q_2u_2\overline{v_2} \dd x &= \langle \partial_\nu u_2, v_2\rangle.
\end{align*}
Setting $v:=v_1-u_2$ and $q_0=q_1-q_2$ we get after subtracting
these identities
\[
\int_\Omega \nabla v\cdot \overline{\nabla v_2} + q_2
v\overline{v_2} + q_0 v_1\overline{v_2} =
\langle(\Lambda_{1}-\Lambda_{2})v_1,v_2 \rangle.
\]
Since $v_2$ solves $(\Delta-q_2)v_2=0$, $v=0$ on $\partial\Omega$
and $v_2=0$ on $\Gamma_0$, we have
\[
\int_\Omega \nabla v \cdot\overline{\nabla v_2} + q_2
v\overline{v_2}=0,
\]
\[
\langle(\Lambda_{1}-\Lambda_{2})v_1,v_2 \rangle=\langle
(\Lambda_{1,\Gamma}-\Lambda_{2,\Gamma})v_1,v_2 \rangle,
\]
and the assertion follows.\qed

In treating inverse boundary value problems, complex geometrical
optics solutions play a very important role. We now describe the complex
geometrical optics solutions we are going to use in our proofs. We will follow the idea in \cite{I07}.
Assume that $q_1$, $q_2\in L^{\infty}(\R^n)$ are compactly supported and are even in $x_n$, i.e.
$$q_1^{\ast}(x_1,\cdots,x_{n-1},x_n)=q_1(x_1,\cdots,x_{n-1},x_n)$$
and
$$q_2^{\ast}(x_1,\cdots,x_{n-1},x_n)=q_2(x_1,\cdots,x_{n-1},x_n).$$ Hereafter, we denote
$$h^{\ast}(x_1,\cdots,x_{n-1},x_n)=h(x_1,\cdots,x_{n-1},-x_n).$$

Given $\xi=(\xi_1,\cdots,\xi_n)\in \R^n$.
Let us first introduce new coordinates obtained by rotating the standard Euclidean coordinates around the $x_n$ axis such that
the representation of $\xi$ in the new coordinates, denoted by $\tilde\xi$, satisfies
$\tilde\xi=(\tilde\xi_1,0,\cdots,0,\tilde\xi_n)$ with $\tilde\xi_1=\sqrt{\xi_1^2+\cdots+\xi^2_{n-1}}$ and $\tilde\xi_n=\xi_n$.
In the following we also denote by $\tilde x$ the representation of $x$ in the new coordinates. Then we define for $\tau>0$
\begin{equation}\label{rho}
\begin{split}
\tilde\rho_1 & :=(\frac{\tilde\xi_1}{2}-\tau\tilde\xi_n,i|\xi|(\frac14+\tau^2)^{1/2},0,\cdots,0,
            \frac{\tilde\xi_n}{2}+\tau\tilde\xi_1),\\
\tilde\rho_2 & :=(\frac{\tilde\xi_1}{2}+\tau\tilde\xi_n,-i|\tilde\xi|(\frac14+\tau^2)^{1/2},0,\cdots,0,
            \frac{\tilde\xi_n}{2}-\tau\tilde\xi_1),
\end{split}
\end{equation}
and let $\rho_1$ and $\rho_2$ be representations of $\tilde\rho_1$ and $\tilde\rho_2$ in the original coordinates. Note that $x_n=\tilde x_n$ and $\sum_{i=1}^n x_iy_i=\sum_{i=1}^n \tilde
x_i\tilde y_i$. It is clear that, for $j=1,2$, $\rho_j\cdot\rho_j=0$ as well as
$\rho^*_j\cdot\rho^*_j=0$ hold.

The construction given in \cite{SU87}
ensures that there are complex geometrical optics solutions
$u_j=e^{i\rho_j\cdot x}(1+w_j)$ of $(\Delta-q_j)u_j=0$ in $\R^n$, $j=1,2$, and the functions $w_j$ satisfy $\|w_j\|_{L^2(K)}\leq C_K \tau^{-1}$ for any compact set $K\subset\R^n$. We then set
\begin{equation}\label{spec_sol}
\begin{split}
v_1(x) &= e^{i\rho_1\cdot x}(1+w_1) - e^{i\rho^*_1\cdot x}(1+w_1^*)\\
v_2(x) &= e^{-i\rho_2\cdot x}(1+w_2) - e^{-i\rho^*_2\cdot x}(1+w_2^*).
\end{split}
\end{equation}
From this definition it is clear that these functions are solutions of
$(\Delta -q_j)v_j=0$ in $\R^n_+$ with $v_j=0$ on $x_n=0$.
\section{Stability estimate for the potential}
Now we are in the position to prove Theorem~\ref{main_bdd}. We first
consider the case (a) where $\Gamma_0$ is a part of a hyperplane. To
construct the special solutions described in the previous section,
we first perform zero extension of $q_1$ and $q_2$ to $\R_n^{+}$ and
then even extension to the whole $\R^n$. As in the last section, we
can construct special geometrical optics solutions $v_j$ of the form
\eqref{spec_sol} to $(\Delta-q_j)v_j=0$ in $\Omega$ for $j=1,2$.
Note that $v_1=v_2=0$ on $\Gamma_0$. We now plug in these solutions
into the identity \eqref{lma_id} and write $q_0=q_1-q_2$. This
gives
\begin{align}\label{eqn:ident}
\begin{split}
&\langle (\Lambda_{1,\Gamma}-\Lambda_{2,\Gamma}) v_1,
        v_2\rangle\\
&= \int_{\Omega} q_0 v_1\overline{v_2}\dd x\\
&= \int_{\Omega} q_0(x)\Big{(}e^{i(\rho_1+\rho_2)\cdot x}(1+w_1)(1+\overline{w_2})+e^{i(\rho^{\ast}_1+\rho^{\ast}_2)\cdot x}(1+w^{\ast}_1)(1+\overline{w^{\ast}_2})\\
&\qquad\quad -e^{i(\rho_1+\rho_2^{\ast})\cdot x}(1+w_1)(1+\overline{w_2^{\ast}})-e^{i(\rho_1^{\ast}+\rho_2)\cdot x}(1+w_1^{\ast})(1+\overline{w_2})\Big{)}\dd x\\
&= \int_{\Omega} q_0(x)(e^{i\xi\cdot x}+e^{i\xi^{\ast}\cdot x})\dd x +\int_{\Omega}q_0(x)f(x,w_1, w_2,
        w_1^*, w_2^*)\dd x\\
&\qquad\quad -\int_{\Omega}q_0(x)\big{(}e^{i(\rho_1+\rho_2^{\ast})\cdot x}
        +e^{i(\rho_1^{\ast}+\rho_2)\cdot x}\big{)}\dd x,
\end{split}
\end{align}
where
\begin{align*}
f\ =&\ e^{i\xi\cdot x}(w_1+\overline{w_2}+w_1\overline{w_2})+ e^{i\xi^*\cdot x}(w^*_1+\overline{w^*_2}+w^*_1\overline{w^*_2})\\
    &-e^{i(\rho_1^{\ast}+\rho_2)\cdot x}(w^*_1+\overline{w_2}+w^*_1\overline{w_2})-e^{i(\rho_1+\rho_2^{\ast})\cdot x}(w_1+\overline{w^*_2}+w_1\overline{w^*_2}).
\end{align*}

The first term on the right hand side of \eqref{eqn:ident} is equal to
\[
\int_{\R^n}q_0(x)e^{i\xi\cdot x}\dd x=\cF q_0(\xi)
\]
because $q_0$ is even in $x_n$. For the second term, we use the estimate
\[
\|w_1\|_2 +\|w^*_1\|_2 +\|\overline{w_2}\|_2 +\| \overline{w^*_2}\|_2 \leq
C\tau^{-1}
\]
to obtain
\begin{equation}\label{eqn:f_est}
\big{|}\int_\Omega q_0 f(x,w_1, w_2,w_1^*, w_2^*)\dd x\big{|}\leq C\|q_0\|_2
\tau^{-1}.
\end{equation}
As for the last term on the right hand side of \eqref{eqn:ident}, we first observe that
\[(\rho_1+\rho_2^{\ast})\cdot x=(\tilde\rho_1+\tilde\rho_2^{\ast})\cdot\tilde x=\tilde\xi_{1}
        \tilde x_{1}+2\tau\tilde\xi_{1}\tilde x_n=\xi'\cdot x'+2\tau|\xi'|x_n
\]
and
\[
(\rho_1^{\ast}+\rho_2)\cdot x=(\tilde\rho_1+\tilde\rho_2^{\ast})\cdot\tilde x=\tilde\xi_{1}
        \tilde x_{1}+2\tau\tilde\xi_{1}\tilde x_n=\xi'\cdot x'-2\tau|\xi'|x_n,
\]
where $\xi'=(\xi_1,\cdots,\xi_{n-1})$ and $x'=(x_1,\cdots,x_{n-1})$.
Therefore, we can write
\[
\int_{\Omega} q_0(x)e^{i(\rho_1+\rho_2^{\ast})\cdot x}\dd x = \cF q_0(\xi',2\tau|\xi'|)
\]
as well as
\[
\int_{\Omega} q_0(x)e^{i(\rho_1^{\ast}+\rho_2)\cdot x}\dd x= \cF q_0(\xi',-2\tau|\xi'|).
\]
The Sobolev embedding and the assumptions on $q_j$ ensure that
$q_0\in C^\alpha(\overline{\Omega})$ for $\alpha=s-\frac{n}{2}$ and
therefore $q_0$ satisfies the assumption of
Lemma~\ref{lma:holder_est}. Applying Lemma~\ref{fourier_est} to
$q_0$ yields that for $\eps<\eps_0$
\begin{equation}\label{eqn:q0est}
|\cF q_0(\xi',2\tau|\xi'|)|+
|\cF q_0(\xi',-2\tau|\xi'|)|
    \leq C(\exp(-\pi\eps^2(1+4\tau^2)|\xi'|^2)+\varepsilon^{\alpha}).
\end{equation}
Finally, we estimate the boundary integral
\begin{align}\label{bi}
\begin{split}
&\left|\int_{\Gamma} (\Lambda_{1,\Gamma}-\Lambda_{2,\Gamma}) v_1
\cdot v_2 \dd \sigma\right|\\
&\le\|\Lambda_{1,\Gamma}-\Lambda_{2,\Gamma}\|_*
        \|v_1\|_{H^{\frac{1}{2}}(\Gamma)}
        \|v_2\|_{H^{\frac{1}{2}}(\Gamma)}\\
&\le\|\Lambda_{1,\Gamma}-\Lambda_{2,\Gamma}\|_*
        \|v_1\|_{H^1(\Omega)}
        \|v_2\|_{H^1(\Omega)}\\
&\le C\exp(|\xi|\tau)\|\Lambda_1-\Lambda_2\|_*.
\end{split}
\end{align}
Combining \eqref{eqn:ident}, \eqref{eqn:f_est}, \eqref{eqn:q0est},
and \eqref{bi} leads to the inequality
\begin{equation}\label{ineq1}
|\cF q_0(\xi)|\le
C\{\exp(|\xi|\tau)\|\Lambda_1-\Lambda_2\|_*+\exp(-\pi\eps^2(1+4\tau^2)|\xi'|^2)+\varepsilon^{\alpha}+\frac{1}{\tau}\}
\end{equation}
for all $\xi\in\R^n$ and $\eps<\eps_0$, where $C$ only depends on a
priori data on the potentials.

Next we would like to estimate the norm of $q_0$ in $H^{-1}$. As
usual, other estimates of $q_0$ in more regular norms can be
obtained by interpolation. To begin, we set $Z_R=\{\xi \in \R^n:\;
|\xi_n|<R \mbox{ and }|\xi'|<R \}$. Note that $B(0,R)\subset
Z_R\subset B(0,cR)$ for some $c>0$. Now we use the a priori
assumption on potentials and \eqref{ineq1} and calculate
\begin{align}\label{ineq2}
\begin{split}
\|q_0\|_{H^{-1}}^2
    \leq &\int_{Z_R}|\cF q_0(\xi)|^2(1+|\xi|^2)^{-1}\dd \xi
        +\int_{{Z_R}^c}|\cF q_0(\xi)|^2(1+|\xi|^2)^{-1}\dd \xi\\
    \leq &\int_{Z_R}|\cF q_0(\xi)|^2(1+|\xi|^2)^{-1}\dd \xi+CR^{-2}\\
    \leq\strut &C\{R^n\exp(cR\tau)\|\Lambda_1-\Lambda_2\|_*^2+R^n\eps^{2\alpha}+R^n\tau^{-2}+R^{-2}\\
    &   +\int_{-R}^{R}\int_{B'(0,R)}\exp(-2\pi\eps^2(1+4\tau^2)|\xi'|^2)\dd
        \xi'\dd \xi_n \},\\
\end{split}
\end{align}
here $B'(x',R)$ denotes the ball in $\R^{n-1}$ with center $x'$ and
radius $R>0$. For the second term on the right hand side of
\eqref{ineq2}, we choose $\eps=(1+4\tau^2)^{-1/4}$ with
$\tau\ge\tau_0\gg 1$ and integrate
\begin{align}\label{ineq3}
\begin{split}
\lefteqn{\int_{-R}^{R}\int_{B'(0,R)}\exp(-2\pi\eps^2(1+4\tau^2)|\xi'|^2)\dd
        \xi'\dd \xi_n}\\
&= 2R\int_{B'(0,R)}
    \exp(-2\pi(1+4\tau^2)^{1/2}|\xi'|^2)\dd \xi'\\
&= 2R\int_{S^{n-2}}\int_{0}^R r^{n-2}
    \exp(-2\pi((1+4\tau^2)^{1/4}r)^2)\dd r\dd \omega\\
&\leq CR (1+4\tau^2)^{-(n-1)/4}\int_{0}^\infty u^{n-2}
    \exp(-2\pi u^2)\dd u\\
&\leq CR\tau^{-(n-1)/2}.
\end{split}
\end{align}
Plugging  \eqref{ineq3} into \eqref{ineq2} with the choice of
$\eps=(1+4\tau^2)^{-1/4}$ we get for $R>1$
\begin{align}\label{ineq5}
\begin{split}
\|q_0\|^2_{H^{-1}}&\leq C\{R^n
\exp(cR\tau)\|\Lambda_1-\Lambda_2\|_*^2+R^n\tau^{-\alpha}+R\tau^{-(n-1)/2}+R^{-2}\}\\
&\leq C\{R^n
\exp(cR\tau)\|\Lambda_1-\Lambda_2\|_*^2+R^n\tau^{-\tilde\alpha}+R^{-2}\},
\end{split}
\end{align}
where $\tilde\alpha=\min\{\alpha,(n-1)/2\}$.

Observing from \eqref{ineq5}, we now choose $\tau$ such that
$R^n\tau^{-\tilde\alpha}=R^{-2},$ namely,
$\tau=R^{(n+2)/\tilde\alpha}$. Substituting such $\tau$ back to
\eqref{ineq5} yields
\begin{equation}\label{ineq6}
\|q_0\|^2_{H^{-1}}\leq C\{R^n
\exp(cR^{\frac{n+2}{\tilde\alpha}+1})\|\Lambda_1-\Lambda_2\|_*^2+R^{-2}\}.\\
\end{equation}
Finally, we choose a suitable $R$ so that
\[
R^n\exp(cR^{\frac{n+2}{\tilde\alpha}+1})\|\Lambda_1-\Lambda_2\|_*^2=R^{-2},
\]
i.e., $R=\big{|}\log\|\Lambda_1-\Lambda_2\|_*\big{|}^{\gamma}$ for
some $0<\gamma=\gamma(n,\tilde\alpha)$. Thus, we obtain from
\eqref{ineq6} that
\begin{equation}\label{ineq10}
\|q_1-q_2\|_{H^{-1}(\Omega)}\leq
C\big{|}\log\|\Lambda_1-\Lambda_2\|_*\big{|}^{-\gamma}.
\end{equation}
The derivation of \eqref{ineq10} is legitimate under the assumption
that $\tau$ is large. To make sure that it is true, we need to take
$R$ sufficiently large, i.e. $R>R_0$ for some large $R_0$.
Consequently, there exists $\tilde\delta>0$ such that if
$\|\Lambda_1-\Lambda_2\|_*<\tilde\delta$ then \eqref{ineq10} holds.
For $\|\Lambda_1-\Lambda_2\|_*\ge\tilde\delta$, \eqref{ineq10} is
automatically true with a suitable constant $C$ when we take into
account the a priori bound \eqref{apr}.

The estimate \eqref{est} is now an easy consequence of the
interpolation theorem. Precisely, let $\epsilon>0$ such that
$s=\frac{n}{2}+2\epsilon$. Using that
$[H^{t_0}(\Omega),H^{t_1}(\Omega)]_\beta = H^t(\Omega)$ with
$t=(1-\beta)t_0+\beta t_1$ (see e.g. \cite[Theorem 1 in
4.3.1]{Tri95}) and the Sobolev embedding theorem, we get
$\|q_1-q_2\|_{L^\infty} \leq
C\|q_1-q_2\|_{H^{\frac{n}{2}+\epsilon}}\leq
C\|q_1-q_2\|^{(1-\beta)}_{H^{t_0}} \|q_1-q_2\|^{\beta}_{H^{t_1}}$.
Setting $t_0=-1$ and $t_1=s$ we end up with
\[
\|q_1-q_2\|_{L^\infty(\Omega)} \leq C
\|q_1-q_2\|^{\frac{s+1-\epsilon}{s+1}}_{H^{-1}(\Omega)}
\]
which yields the desired estimate \eqref{est} with
$\sigma=\gamma\frac{s+1-\epsilon}{s+1}$.

We now turn to case (b). With a suitable translation and rotation,
it suffices to assume $a=(0,\cdots,0,R)$ and
$0\notin\overline\Omega$. As in \cite{I07}, we shall use Kelvin's
transform:
\begin{equation}\label{kel}
y=\left(\frac{2R}{|x|}\right )^2x\quad\text{and}\quad
x=\left(\frac{2R}{|y|}\right)^2y.
\end{equation}
Let $$\tilde u(y)=\left(\frac{2R}{|y|}\right)^{n-2}u(x(y)),$$ then
$$\left(\frac{|y|}{2R}\right)^{n+2}\Delta_y\tilde u(y)=\Delta_xu(x).$$
Denote by $\tilde\Omega$ the transformed domain of $\Omega$. In view of
this transform, $\Gamma_0$ now becomes $\tilde\Gamma_0\subset
\{y_n=2R\}$ and $\Gamma$ is transformed to $\tilde\Gamma$ and
$\tilde\Gamma=\partial\tilde\Omega\cap\{y_n>2R\}$. On the other
hand, if $u(x)$ satisfies $\Delta u-q(x)u=0$ in $\Omega$, then
$\tilde u$ satisfies
\begin{equation}\label{teq}
\Delta \tilde u-\tilde q\tilde u=0\quad\text{in}\quad\tilde\Omega,
\end{equation}
where
\[
\tilde q(y)=\left(\frac{2R}{|y|}\right)^{4}q(x(y)).
\]
Therefore,
for \eqref{teq} we can define the partial Dirichlet-to-Neumann map
$\tilde\Lambda_{\tilde q,\tilde\Gamma}$ acting boundary functions
with homogeneous data on $\tilde\Gamma_0$.

We now want to find the relation between $\Lambda_{q,\Gamma}$ and
$\tilde\Lambda_{\tilde q,\tilde\Gamma}$. It is easy to see that for
$f,g\in H^{1/2}_0(\Gamma)$
\[
\langle\Lambda_{q,\Gamma} f,g\rangle=\int_{\Omega}(\nabla
u\cdot\nabla\overline{v}+qu\overline{v})\dd x,
\]
where $u$ solves
\begin{align*}
\Delta u-qu&=0\quad\text{in}\quad\Omega,\\
u&=f\quad\text{on}\quad\partial\Gamma
\end{align*}
and $v\in H^1(\Omega)$ with $v|_{\partial\Omega}=g$. Defining
\[
\tilde f=\left(\frac{2R}{|y|}\right)^{n-2}\Big{|}_{\partial\tilde\Omega}f,\quad
\tilde g=\left(\frac{2R}{|y|}\right)^{n-2}\Big{|}_{\partial\tilde\Omega}g,
\]
and
\[
\quad
\tilde v(y)=\left(\frac{2R}{|y|}\right)^{n-2}v(x(y)).
\] 
Then we have
$\tilde f,\tilde g\in H^{1/2}_0(\tilde\Gamma)$ and
\begin{equation*}
\langle\Lambda_{q,\Gamma} f,g\rangle=\langle\tilde\Lambda_{\tilde
q,\tilde\Gamma}\tilde f,\tilde g\rangle,
\end{equation*}
in particular,
\begin{equation}\label{2lam}
\langle(\Lambda_{q_1,\Gamma}-\Lambda_{q_2,\Gamma})
f,g\rangle=\langle(\tilde\Lambda_{\tilde
q_1,\tilde\Gamma}-\tilde\Lambda_{\tilde q_2,\tilde\Gamma})\tilde
f,\tilde g\rangle.
\end{equation}
With the assumption $0\notin\overline{\Omega}$, the change of
coordinates $x\to y$ by \eqref{kel} is a diffeomorphism from
$\overline{\Omega}$ onto $\overline{\tilde\Omega}$. Note that
$(2R/|y|)^{n-2}$ is a positive smooth function on
$\partial\tilde\Omega$. Recall a fundamental fact from Functional
Analysis:
\begin{equation}\label{norm}
\|\Lambda_{q_1,\Gamma}-\Lambda_{q_2,\Gamma}\|_{\ast}=\sup\left\{\frac{|\langle(\Lambda_{q_1,\Gamma}-\Lambda_{q_2,\Gamma})
f,g\rangle|}{\|f\|_{H^{1/2}_0(\Gamma)}\|g\|_{H^{1/2}_0(\Gamma)}}\ :\
f,g\in H^{1/2}_0(\Gamma)\right\}.
\end{equation}
The same formula holds for $\|\tilde\Lambda_{\tilde
q_1,\tilde\Gamma}-\tilde\Lambda_{\tilde q_2,\tilde\Gamma}\|_{\ast}$.
On the other hand, it is not difficult to check that
$\|f\|_{H^{1/2}_0(\Gamma)}$ and $\|\tilde
f\|_{H^{1/2}_0(\tilde\Gamma)}$, $\|g\|_{H^{1/2}_0(\Gamma)}$ and
$\|\tilde g\|_{H^{1/2}_0(\tilde\Gamma)}$ are equivalent, namely,
there exists $C$ depending on $\partial\Omega$ such that
\begin{align}\label{equiv}
\begin{split}
&\frac 1C \|f\|_{H^{1/2}_0(\Gamma)}\le\|\tilde
f\|_{H^{1/2}_0(\tilde\Gamma)}\le C\|f\|_{H^{1/2}_0(\Gamma)},\\
&\frac 1C \|g\|_{H^{1/2}_0(\Gamma)}\le\|\tilde
g\|_{H^{1/2}_0(\tilde\Gamma)}\le C\|g\|_{H^{1/2}_0(\Gamma)}.
\end{split}
\end{align}
Putting together \eqref{2lam}, \eqref{norm}, and \eqref{equiv} leads
to
\begin{equation}\label{lest}
\|\tilde\Lambda_{\tilde q_1,\tilde\Gamma}-\tilde\Lambda_{\tilde
q_2,\tilde\Gamma}\|_{\ast}\le
C\|\Lambda_{q_1,\Gamma}-\Lambda_{q_2,\Gamma}\|_{\ast}
\end{equation}
with $C$ only depending on $\partial\Omega$.

With all the preparations described above, we use case (a) for
the domain $\tilde\Omega$ with the partial Dirichlet-to-Neumann map
$\tilde\Lambda_{\tilde q,\tilde\Gamma}$. Therefore, we immediately
obtain the estimate:
\[
\|\tilde q_1-\tilde q_2\|_{L^\infty(\tilde\Omega)} \leq C
\big{|}\log \|\tilde\Lambda_{\tilde q_1,\tilde\Gamma}-
    \tilde\Lambda_{\tilde q_2,\tilde\Gamma}\|_{*} \big{|}^{-\sigma}.
\]
Finally, rewinding $\tilde q$ and using \eqref{lest} yields the
estimate \eqref{est}.

\section{Stability estimate for the conductivity}

We aim to prove Corollary~\ref{cor1} in this section. We recall the
following well-known relation: let
$q=\frac{\Delta\sqrt{\gamma}}{\sqrt{\gamma}}$ then
\[
\Lambda_{q,\Gamma}(f)=\gamma^{-1/2}|_{\Gamma}\Lambda_{\gamma,\Gamma}(\gamma^{-1/2}|_{\Gamma}f)+\frac{1}{2}(\gamma^{-1}\partial_{\nu}\gamma)|_{\Gamma}f.
\]
In view of the a priori assumption \eqref{sameb}, we have that
\begin{equation*}
(\Lambda_{q_1,\Gamma}-\Lambda_{q_2,\Gamma})(f)=\gamma^{-1/2}|_{\Gamma}(\Lambda_{\gamma_1,\Gamma}-\Lambda_{\gamma_2,\Gamma})(\gamma^{-1/2}|_{\Gamma}f)
\end{equation*}
where
$\gamma^{-1/2}|_{\Gamma}:=\gamma_1^{-1/2}|_{\Gamma}=\gamma_2^{-1/2}|_{\Gamma}$,
which implies
\begin{equation}\label{qgamma}
\|\Lambda_{q_1,\Gamma}-\Lambda_{q_2,\Gamma}\|_*\le
C\|\Lambda_{\gamma_1,\Gamma}-\Lambda_{\gamma_2,\Gamma}\|_*
\end{equation}
for some $C=C(N)>0$. Hereafter, we set
$q_j=\frac{\Delta\sqrt{\gamma_j}}{\sqrt{\gamma_j}}$, $j=1,2$. The
regularity assumption \eqref{ap} and Sobolev's embedding theorem
imply that $q_1,q_2\in C^1(\overline{\Omega})$. Using this and
\eqref{sameb}, we conclude that $\hat q_1-\hat q_2$ satisfies the
assumptions of Lemma~\ref{lma:holder_est} with $\alpha=1$.
Therefore, Theorem~\ref{main_bdd} and
\eqref{qgamma} imply that
\begin{equation}\label{mdest}
\|q_1-q_2\|_{L^\infty(\Omega)} \leq C \big{|}\log
\|\Lambda_{\gamma_1,\Gamma}-
    \Lambda_{\gamma_2,\Gamma}\|_{*} \big{|}^{-\sigma_1}
\end{equation}
where $C$ depend on $\Omega, N, n, s$ and $\sigma_1$ depend on
$n,s$. Next, we recall from \cite[(26) on page 168]{A88} that
\begin{equation}\label{qg}
\|\gamma_1-\gamma_2\|_{L^{\infty}(\Omega)}\le
C\|q_1-q_2\|^{\sigma_2}_{L^{\infty}(\Omega)}
\end{equation}
for some $0<\sigma_2<1$, where $C=C(N,\Omega)$ and
$\sigma_2=\sigma_2(n,s)$. Finally, putting together \eqref{mdest}
and \eqref{qg} yields \eqref{est111} with $\sigma=\sigma_1\sigma_2$
and the proof of Corollary~\ref{cor1} is complete.


\vspace{.3cm}\noindent
\begin{thebibliography}{999999}
\bibitem[Al88]{A88}
G. Alessandrini, \emph{Stable determination of conductivity by
boundary measurements}, Appl. Anal. \textbf{27} (1988), no.~1-3,
153-172.

\bibitem[Al90]{Al90}
G. Alessandrini, \emph{Singular solutions of elliptic equations and 
the determination of conductivity by boundary measurements}, 
J. Differential Equations, \textbf{84} (1990), 252-272.

\bibitem[AP06]{AP06}
K. Astala and L. P\"aiv\"arinta, \emph{Calderon's inverse
conductivity problem in the plane}, Ann. of Math. (2), \textbf{163}
(2006), 265-299.

\bibitem[BU02]{BU02}
A. L. Bukhgeim and G. Uhlmann, \emph{Recovering a potential from
partial Cauchy data}, Comm. Partial Differential Equations
\textbf{27} (2002), 653-668.

\bibitem[Ca80]{C80}
A. Calder\'on, \emph{On an inverse boundary value problem}, Seminar
on Numerical Analysis and its Applications to Continuum Physics,
Soc. Brasileira de Matem\'{a}tica, R\'io de Janeiro (1980), 65-73.

\bibitem[FKSU07]{FKSU07}
D. Dos Santos Ferreira, C.E. Kenig, J. Sj\"{o}strand, and G.
Uhlmann, \emph{Determining the magnetic Schr\"{o}dinger operator
from partial Cauchy data}, Comm. Math. Phys., to appear.

\bibitem[GU01]{GU01}
A.~Greenleaf and G.~Uhlmann, \emph{Local uniqueness for the
  {D}irichlet-to-{N}eumann map via the two-plane transform}, Duke Math. J.
  \textbf{108} (2001), 599-617.

\bibitem[HW06]{HW06}
H. Heck and J.-N. Wang, \emph{Stability estimates for the inverse
boundary value problem by partial Cauchy data}, Inverse Problems
\textbf{22} (2006), 1787-1796.

\bibitem[IU04]{IU04}
H. Isozaki and G. Uhlmann, \emph{Hyperbolic geometry and the local
Dirichlet-to-Neumann map}, Advances in Math., {\bf 188} (2004),
294-314.

\bibitem[Is07]{I07}
V. Isakov, \emph{On uniqueness in the inverse conductivity problem
with local data}, Inverse Problems and Imaging, \textbf{1} (2007),
95-105.

\bibitem[KSU05]{KSU05}
C.E. Kenig, J. Sj\"{o}strand, and G. Uhlmann, \emph{The Calder\'{o}n
problem with partial data}, to appear in Ann. of Math.

\bibitem[Ma01]{M01}
N. Mandache,  \emph{Exponential instability in an inverse problem
for the Schrödinger equation}, Inverse Problems \textbf{17}  (2001),
no. 5, 1435-1444.

\bibitem[Na88]{N88}
A. Nachman, \emph{Reconstructions from boundary measurements},  Ann.
of Math. (2) \textbf{128} (1988),  no. 3, 531-576.

\bibitem[Na96]{N96}
\bysame, \emph{Global uniqueness for a two-dimensional inverse
boundary value problem}, Ann. of Math. (2), \textbf{143} (1996),
71-96.

\bibitem[SU87]{SU87}
J.~Sylvester and G.~Uhlmann, \emph{A global uniqueness theorem for
an inverse  boundary value problem}, Ann. of Math. (2) \textbf{125}
(1987), 153-169.

\bibitem[SU88]{SU88}
J. Sylvester and G. Uhlmann, \emph{Inverse boundary value problems at 
the boundary---continuous dependence}, 
Comm. Pure Appl. Math., \textbf{41} (1988), 197-219.

\bibitem[Tr95]{Tri95}
H.~Triebel, \emph{{Interpolation theory, function spaces,
differential
  operators}}, 2nd ed., {Johann Ambrosius Barth}, Heidelberg, Leipzig, 1995.

\bibitem[Ve99]{Ve99}
S.~Vessella, \emph{A continuous dependence result in the analytic
continuation problem}, Forum Math. \textbf{11} (1999), 695-703.
\end{thebibliography}
\end{document}